\documentclass[12pt]{article}
\usepackage{tikz}
\usetikzlibrary{tikzmark}
\usepackage{tikz-opm}
 \usepackage{amsmath, amsthm, amscd, amsfonts,graphicx}
\usepackage[matrix,arrow]{xy}

\newtheorem{thm}{Theorem}[section]
\newtheorem{dfn}{Definition}[section]
\newtheorem{lemma}{Lemma}[section]
\newtheorem{prop}{Proposition}[section]
\newtheorem{rem}{Remark}[section]
\newtheorem{ex}{Example}[section]
\newtheorem{prf}{Proof}[section]

\newcommand{\blem}{\begin{lemma}}
\newcommand{\elem}{\end{lemma}}

\newcommand{\bprf}{\begin{prf}}
\newcommand{\eprf}{\end{prf}}
\newcommand{\bthm}{\begin{thm}}
\newcommand{\ethm}{\end{thm}}

\newcommand{\bdf}{\begin{dfn}\rm}
\newcommand{\edf}{\end{dfn}}

\newcommand{\bpro}{\begin{prop}}
\newcommand{\epro}{\end{prop}}

\newcommand{\brem}{\begin{rem}\rm}
\newcommand{\erem}{\end{rem}}

\newcommand{\bex}{\begin{ex}\rm}
\newcommand{\eex}{\end{ex}}

\newcommand{\zerothemall}{
\setcounter{prop}{0} \setcounter{lemma}{0} \setcounter{rem}{0}
\setcounter{thm}{0} \setcounter{ex}{0} \setcounter{dfn}{0}
\setcounter{equation}{0} }

\title{An Analysis of A Fishing Model with Nonlinear Harvesting Function}
\author{E. SOROURI, M. ESHAGHI GORDJI, R. MEMARBASHI \\
Department of Mathematics, Semnan University,\\ P. O. Box 35195-363,
Semnan, Iran.\\
sorouri.e@semnan.ac.ir,\\
 meshaghi@semnan.ac.ir, \\
 r\_memarbashi@semnan.ac.ir }
\date{}

\begin{document}
\maketitle

\begin{abstract}
\noindent In this study, considering the importance of how to exploit renewable natural resources, we analyze a fishing model with nonlinear harvesting function in which the players at the equilibrium point do a  static game 
 with complete information  that, according to the calculations, will cause a waste of energy for both players
and so the selection of cooperative strategies along with the agreement between the players is the result of this research.

\vspace{5mm}

\noindent{\it Keywords: natural resource management, game theory, bioeconomic models, nonlinear dynamical systems, fishery, harvest function\\
JLE classification:Q22, Q57\\
AMS classification:91B76, 91A80}

\end{abstract}
\section{Introduction}
\zerothemall 
Since game theory examines situations in which decision-makers interact,  this theory has many applications in the commercial competition between individuals, companies, and countries (see \cite{Bernheim},\cite{Dixit},\cite{Gibbons},\cite{Osborne} and \cite{OsborneM}).\\
For example, using this theory, we can examine how to use renewable natural resources and different strategies. 
Consider a river, lake or sea exploited by fishermen or companies.
 If the number of fishermen is increased or more fish are harvested,  it will lead to the
extinction of the generation of fishes in that source.

If fishering from the source is done  only by a fisherman, he(or she) will consider the amount of current harvest because
 he(or she) knows that the amount of fish may be reduced and harvesting in the future will require more effort and cost, so 
he(or she) considers the effects of fishing today in the effort and cost of future fishing.
Now if two individuals or companies fish from a natural resource, then fishing more today that leads to  more cost and effort in the future and it will be shared between them.
Therefore, there is a motivation for more harvesting in the present  by either fisherman.\\
The importance of exploiting renewable natural resources and paying attention to the above points have  led to a widespread examination of the issue of fishing and harvest management strategies that prevent the extinction of species (see \cite{Agnew},\cite{Bisch},\cite{Clark},\cite{Kaitala },\cite{VKaitala},\cite{KaitalaM} and \cite{Sumaila}).\\
Of course, a large number of these studies have been investigated in Dynamic Systems (see \cite{ClarkC}, \cite{Cohen}, \cite{Dubey}, \cite{Schaefer}, \cite{MSchaefer} and \cite{Van}).\\
The logistic equation with density-dependent harvesting   (see \cite{Cooke}) is 
$$\dfrac{dN}{dt}=rN(1-\frac{N}{K})-H(t).$$
where $N$ is the population biomass of fish at time $t$, $r$ is the intrinsic rate of growth of the population, $K$ is the carrying  capacity, and $H(t)$ as the harvest function is 
$$H(t)=qEN(t)$$
where $E$ is the fishing effort, the intensity  of the human activities to extract the fish and  
 $q\geq0$ is the catchability coefficient which is defined as the fraction of the population fished by a unit of effort.\\
In the above fishing mdel and many other models, the harvest function is linear in terms of $N(t)$ but we want to examine the effect of the  nonlinear harvest  function of $H(t)=qE(N(t))^{2}$
and since that  fishermen or companies usually harvest  individually and based on their own profits, at the point of equilibrium of the system we consider a static game with complete information and  as a result  we calculate the amount of the waste of effort  and energy.

\section{Model formulation and basic properties}
\zerothemall

In this model, we consider the relation between net growth, $W$, and the carrying  capacity, $K$,  a logistic growth function. So,  we have
\begin{equation}\label{2.1}
W=rN(1-\frac{N}{K}).
\end{equation}

Assuming a nonlinear harvest function, $H=qEN^{2}$, the dynamics of this model is

$$\dfrac{dN}{dt}=rN(1-\frac{N}{K})-qEN^{2}.$$
It can be written as follows

$$\dfrac{dN}{dt}=rN(1-\frac{N}{K_{0}})$$
where $K_{0}=\frac{rK}{r+qEK}$. This system has a trivial equilibrium, $N=0$,  and an non-trivial equilibrium, $N=K_{0}=\frac{rK}{r+qEK}$.\\
Considering $f(N)=rN(1-\frac{N}{K})-qEN^{2}$, we have $\dfrac{dN}{dt}=f(N)$. Since $\dfrac{df(N)}{dN}=r-2\frac{rN}{K}-2qEN$, then  $\dfrac{df(0)}{dN}=r>0$ so  the trivial equilibrium point is unstable, and at the non-trivial  equilibrium $\dfrac{df(K_{0})}{dN}=\frac{-r^{2}-qEK}{r+qEK}<0$ which shows the stability of this point.

On the other hand the differential equation of this system is a separable type and it is solved by a method called $ separation  \  of     \ variables$.\\
By solving this equation and taking $B=\frac{N(0)}{\vert  K_{0}-N(0)\vert}$, it follows that if $K_{0}-N>0$ then $N(t)=\frac{BK_{0} e^{rt}}{1+B e^{rt}}$  is  the solution of  this differential equation
and since $\lim_{t\rightarrow+\infty} N(t)=K_{0}$, it implies that the non-trivial  equilibrium is asymptotically stable  and  if $K_{0}-N<0$ then $N(t)=\frac{-BK_{0} e^{rt}}{1-B e^{rt}}$ is  the solution of  the differential equation
and $\lim_{t\rightarrow+\infty} N(t)=K_{0}$  implies  that the non-trivial  equilibrium is asymptotically stable.

According to the above, when the system reaches the  equilibrium point, where $H=W$ (see Figure 1), we have 
\begin{equation}\label{2.2}
N=\frac{rK}{r+qEK}.
\end{equation}

Which shows the relation between the equilibrium mass and  the effort. By substituting (\ref{2.2})  in $H(t)$,   the relation between the level of effort and the harvest  is 
\begin{equation}\label{2.3}
H=qE(\frac{rK}{(r+qEK)})^{2}.
\end{equation}
We assume that there are two fishing companies, $A_{1}$ and $A_{2}$,  which use this source separately. 
They as players  do a static game with complete information in devoting the amount of effort to harvest.
If $E_{1}$ and $E_{2}$ respectively represent the level of effort of players $A_{1}$ and $A_{2}$ then
 the total effort to harvest  from this source is $E_{T}=E_{1}+E_{2}$  and the total harvest of this effort is
$H_{T}=qE_{T}(\frac{rK}{(r+qE_{T}K)})^{2}$.\\
In this model, we assume that the share of each company (player) is equal to  its share of total effort in other words
$$H_{1}=\frac{E_{1}}{E_{T}}H_{T}=\frac{E_{1}}{E_{T}}qE_{T}(\frac{rK}{(r+qE_{T}K)})^{2}=qE_{1}(\frac{rK}{(r+q(E_{1}+E_{2})K)})^{2},$$

$$H_{2}=\frac{E_{2}}{E_{T}}H_{T}=\frac{E_{2}}{E_{T}}qE_{T}(\frac{rK}{(r+qE_{T}K)})^{2}=qE_{2}(\frac{rK}{(r+q(E_{1}+E_{2})K)})^{2}$$

where $H_{T}=H_{1}+H_{2}$.\\
If each unit of the harvest in the market has  a value of $P$ and one unit of effort has a cost $C$, then the outcome of the players is
$$U_{A_{1}}(E_{1},E_{2})=PH_{1}-CE_{1}=\frac{qPE_{1}r^{2}K^{2}}{(r+q(E_{1}+E_{2})K)^{2}}-CE_{1}, $$
$$U_{A_{2}}(E_{1},E_{2})=PH_{2}-CE_{2}=\frac{qPE_{2}r^{2}K^{2}}{(r+q(E_{1}+E_{2})K)^{2}}-CE_{2}. $$
To find the Nash equilibrium we use  the best response method.  In this way, we first do the calculations for player 1. So, assuming that $E_{2}=\overline{E_{2}}$ is positive and  constant, we have
$$U_{A_{1}}(E_{1},\overline{E_{2}})=0\Longrightarrow E_{1}=0 \vee  E_{1}=\sqrt{\frac{Pr^{2}}{qC}}-\beta.$$
where $\beta=(\frac{r}{qK}+\overline{E_{2}})$. But 
\begin{equation}\label{2.4}
 E_{1}=\sqrt{\frac{Pr^{2}}{qC}}-\beta>0 \Longleftrightarrow \overline{E_{2}}< r(\sqrt{\frac{P}{qC}}-\frac{1}{qK}).
\end{equation}
On the other hand, if 
$\dfrac{\partial U_{A_{1}}(E_{1},\overline{E_{2}})}{\partial E_{1}}=\frac{qr^{2}PK^{2}(r+q(E_{1}+\overline{E_{2}})K)-2r^{2}q^{2}PK^{3}E_{1}}{(r+q(E_{1}+\overline{E_{2}})K)^{3}}-c=0$
then  by a simple calculations we have
\begin{equation}\label{2.5}
E_{1}^{3}+3E_{1}^{2}\beta+(3\beta^{2}	+\frac{r^{2}P}{qC})E_{1}=\frac{r^{2}P}{qC}\beta-\beta^{3}
\end{equation}
where $\beta=(\frac{r}{qK}+\overline{E_{2}})$.

We consider  $f_{1}(E_{1})=E_{1}^{3}+3E_{1}^{2}\beta+(3\beta^{2}	+\frac{r^{2}P}{qC})E_{1}$ and $f_{2}(E_{1})=\frac{r^{2}P}{qC}\beta-\beta^{3}$ that  according to (\ref{2.4})  and $\beta=\frac{r}{qK}+\overline{E_{2}}>0$  always  $f_{2}(E_{1})=\frac{r^{2}P}{qC}\beta-\beta^{3}>0$.\\
On the other hand, $\lim_{E_{1}\rightarrow+\infty}f_{1}(E_{1})=+\infty$ and  $\lim_{E_{1}\rightarrow-\infty}f_{1}(E_{1})=-\infty$
and also   $f_{1}(E_{1})=E_{1}(E_{1}^{2}+3\beta E_{1}+(3\beta^{2}	+\frac{r^{2}P}{qC}))$ that $E_{1}^{2}+3\beta E_{1}+(3\beta^{2}	+\frac{r^{2}P}{qC})$ has $\Delta=-3\beta^{2}-4\frac{r^{2}P}{qC}<0$
then $f_{1}(E_{1})$ only has a real root $E_{1}=0$.\\ 
Since  $\dfrac{df_{1}(E_{1})}{dE_{1}}=3E_{1}^{2}+6 \beta E_{1}+(3\beta^{2}+\frac{r^{2}P}{qC})$
then $f_{1}(E_{1})$ for $E_{1}>0$ is always an increasing function.\\
According to the above, functions $f_{1}$ and $f_{2}$ for $E_{1}>0$ intersect each other exactly at one point  then 
$U_{A_{1}}(E_{1},\overline{E_{2}})$ for $0<\overline{E_{2}}<r(\sqrt{\frac{P}{qC}}-\frac{1}{qK})$ has   exactly a local extremum that we  show it with $E_{1}^{\circ}$.

In order to determine the type of this extremum, we use the first derivative test for  $U_{A_{1}}(E_{1},\overline{E_{2}})$ in the neighborhood of this point, $E_{1}^{\circ}$.
Since $E_{1}^{\circ}>0$ , we consider $\epsilon>0$ such  that $E_{1}^{\circ}-\epsilon>0$ in this case
\begin{flushleft}

$\frac{d U_{A_{1}}(E_{1}^{\circ},\overline{E_{2}})}{dE_{1}}=
\frac{-C(qKE_{1}^{\circ}+(r+qK\overline{E_{2}}))^{3}-q^{2}r^{2}PK^{3}E_{1}^{\circ}+qr^{2}PK^{2}(r+qK\overline{E_{2}})}{(qKE_{1}^{\circ}+(r+qK\overline{E_{2}}))^{3}}=0$,

$\frac{d U_{A_{1}}(E_{1}^{\circ}+\epsilon,\overline{E_{2}})}{dE_{1}}=
\frac{-C(qK(E_{1}^{\circ}+\epsilon)+(r+qK\overline{E_{2}}))^{3}-q^{2}r^{2}PK^{3}(E_{1}^{\circ}+\epsilon)+qr^{2}PK^{2}(r+qK\overline{E_{2}})}{(qK(E_{1}^{\circ}+\epsilon)+(r+qK\overline{E_{2}}))^{3}}$
$=\frac{-Cq^{3}K^{3}\epsilon^{3}-3qKC\theta_{1}\epsilon-3q^{2}K^{2}C\theta_{1}\epsilon^{2}-q^{2}r^{2}K^{3}P\epsilon}{(qK(E_{1}^{\circ}+\epsilon)+(r+qK\overline{E_{2}}))^{3}}<0$
\end{flushleft}
where $\theta_{1}=qKE_{1}^{\circ}+(r+qK\overline{E_{2}})>0$,
\begin{flushleft}
$\frac{d U_{A_{1}}(E_{1}^{\circ}-\epsilon,\overline{E_{2}})}{d E_{1}}=\frac{Cq^{3}K^{3}\epsilon^{3}+3qKC\theta_{1}\epsilon(qK(E_{1}^{\circ}-\epsilon)
+(r+qK\overline{E_{2}}))+q^{2}r^{2}K^{3}P\epsilon}{(qK(E_{1}^{\circ}-\epsilon)+(r+qK\overline{E_{2}}))^{3}}>0$.
\end{flushleft}

Therefore, according to the first derivative test, $E_{1}^{\circ}$ is 
a relative maximal point for $U_{A_{1}}(E_{1},\overline{E_{2}})$.
According to the above, the best response function of player 1 is
\begin{equation}\label{2.6}
E_{1}=B_{1}(E_{2})= \left\{ \begin{array}{rcl}
E_{1}^{\circ} & \mbox{for}
& E_{2}<r(\sqrt{\frac{P}{qC}}-\frac{1}{qK})\\ 0 & \mbox{for} &  E_{2}\geq r(\sqrt{\frac{P}{qC}}-\frac{1}{qK}). 
\end{array}\right.
\end{equation}

Now we would like to calculate $E_{1}^{\circ}$ as the root  of Equation (\ref{2.5}). For 

$$E_{1}^{3}+3E_{1}^{2}\beta+(3\beta^{2}	+\frac{r^{2}P}{qC})E_{1}+(\beta^{3}-\frac{r^{2}P}{qC}\beta)=0$$ we consider

\begin{center}
$a=3\beta$, $b=3\beta^{2}	+\frac{r^{2}P}{qC}$, and $c=\beta^{3}-\frac{r^{2}P}{qC}\beta$
\end{center}
 then according to the calculation method of the roots  of  a third-order polynomial, we have
\begin{center}
$p=b-\frac{a^{2}}{3}=\frac{r^{2}P}{qC}$, $q=\frac{2a^{3}}{27}-\frac{ab}{3}+c=-2(\frac{r^{2}P\beta}{qC})$ and $\Delta=\frac{q^{2}}{4}+\frac{p^{3}}{27}=\frac{r^{4}P^{2}}{q^{2}C^{2}}(\beta^{2}+\frac{1}{27}\frac{r^{2}P}{qC})>0.$
\end{center}

Since $\Delta>0$, then (\ref{2.5}) has only one real  root  that we already showed it  with $E_{1}^{\circ}$  and this root  is

\begin{flushleft}
$E_{1}^{\circ}=(-\frac{q}{2}+\sqrt{\Delta})^{\frac{1}{3}}+(-\frac{q}{2}-\sqrt{\Delta})^{\frac{1}{3}}-\frac{a}{3}$
$=((\frac{r^{2}P}{qC})(\beta+\sqrt{\beta^{2}+\frac{1}{27}\frac{r^{2}P}{qC}}))^{\frac{1}{3}}+((\frac{r^{2}P}{qC})(\beta-\sqrt{\beta^{2}+\frac{1}{27}\frac{r^{2}P}{qC}}))^{\frac{1}{3}}-\beta.$

\end{flushleft}

Therefore
\begin{flushleft}
$E_{1}=B_{1}(E_{2})$
$ =\left\{ \begin{array}{rcl}
((\frac{r^{2}P}{qC})(\beta+\sqrt{\beta^{2}+\frac{1}{27}\frac{r^{2}P}{qC}}))^{\frac{1}{3}}+((\frac{r^{2}P}{qC})(\beta-\sqrt{\beta^{2}+\frac{1}{27}\frac{r^{2}P}{qC}}))^{\frac{1}{3}}-\beta & \mbox{for}
& E_{2}<r(\sqrt{\frac{P}{qC}}-\frac{1}{qK})\\ 0 & \mbox{for} &  E_{2}\geq r(\sqrt{\frac{P}{qC}}-\frac{1}{qK}) 
\end{array}\right.$
\end{flushleft}
where $\beta=\frac{r}{qK}+E_{2}.$\\
According to the symmetry of the game and with a completely similar discussion for the second player,  we conclude that

\begin{flushleft}
$E_{2}=B_{2}(E_{1})$
$ =\left\{ \begin{array}{rcl}
((\frac{r^{2}P}{qC})(\beta_{*}+\sqrt{\beta_{*}^{2}+\frac{1}{27}\frac{r^{2}P}{qC}}))^{\frac{1}{3}}+((\frac{r^{2}P}{qC})(\beta_{*}-\sqrt{\beta_{*}^{2}+\frac{1}{27}\frac{r^{2}P}{qC}}))^{\frac{1}{3}}-\beta_{*}   & \mbox{for}
& E_{1}<r(\sqrt{\frac{P}{qC}}-\frac{1}{qK})\\ 0 & \mbox{for} &  E_{1}\geq r(\sqrt{\frac{P}{qC}}-\frac{1}{qK}). 
\end{array}\right.$
\end{flushleft}

where $\beta_{*}=\frac{r}{qK}+E_{1}.$
According to the functions $E_{1}$ and $E_{2}$, we can obtain Nash equilibrium for  this model from  solutions of the  following system

\begin{equation}\label{2.7}
\left\{ \begin{array}{rcl}
E_{1}=
((\frac{r^{2}P}{qC})(\beta+\sqrt{\beta^{2}+\frac{1}{27}\frac{r^{2}P}{qC}}))^{\frac{1}{3}}+((\frac{r^{2}P}{qC})(\beta-\sqrt{\beta^{2}+\frac{1}{27}\frac{r^{2}P}{qC}}))^{\frac{1}{3}}-\beta & \mbox{that }
&\beta=\frac{r}{qK}+E_{2}\\
E_{2}=((\frac{r^{2}P}{qC})(\beta_{*}+\sqrt{\beta_{*}^{2}+\frac{1}{27}\frac{r^{2}P}{qC}}))^{\frac{1}{3}}+((\frac{r^{2}P}{qC})(\beta_{*}-\sqrt{\beta_{*}^{2}+\frac{1}{27}\frac{r^{2}P}{qC}}))^{\frac{1}{3}}-\beta_{*}   & \mbox{that} &\beta_{*}=\frac{r}{qK}+E_{1}.

\end{array}\right.
\end{equation}
Since the above system is complicated  based on the parameters  $r$, $q$, $P$, $C$ and $K$, to make a clearer analysis and easier understanding of the problem, we consider the following values for the parameters
$$r=P=1, q=C=\frac{1}{\sqrt{27}} \ and \ K=\sqrt{27}.$$
In this case, the system (\ref{2.7}) is 
\begin{equation}\label{2.8}
\left\{ \begin{array}{r}
E_{1}=3(((1+E_{2})+\sqrt{(1+E_{2})^{2}+1})^{\frac{1}{3}}+((1+E_{2})-\sqrt{(1+E_{2})^{2}+1})^{\frac{1}{3}})-1-E_{2}\\
E_{2}=3(((1+E_{1})+\sqrt{(1+E_{1})^{2}+1})^{\frac{1}{3}}+((1+E_{1})-\sqrt{(1+E_{1})^{2}+1})^{\frac{1}{3}})-1-E_{1}
\end{array}\right.
\end{equation}

We have solved the above nonlinear system with the  iterative method,  it has an approximate  unique solution $(E_{1},E_{2})=(1,1)$, with very little error.
This unique solution is the Nash equilibrium of this model  that  we represent with $(E_{1}^{*},E_{2}^{*})$ therefore $(E_{1}^{*},E_{2}^{*})=(1,1)$. \\
The total effort and harvest  in the Nash equilibrium, $E_{T}$ and $H_{T}$, are
$$E_{T}=E_{1}^{*}+E_{2}^{*}\cong 1+1=2$$ and $$ H_{T}=H_{1}^{*}+H_{2}^{*}$$
where $$H_{1}^{*}=qE_{1}^{*}(\frac{rK}{(r+q(E_{1}^{*}+E_{2}^{*})K)})^{2}=
\frac{\sqrt{27}}{9}=0.5773502691\cong 0.577$$
and similarly $$H_{2}^{*}\cong 0.577$$ then $H_{T}\cong 1.15.$\\
On the other hand, by substituting in (\ref{2.3}) we have 
$$H_{T}=qE_{T}(\frac{rK}{(r+qE_{T}K)})^{2}=\frac{\sqrt{27}E_{T}}{(1+E_{T})^{2}}\cong 1.15.$$

But the recent relation is almost equivalent to 
$$1.15E_{T}^{2}-2.9E_{T}+1.15=0.$$
For this equation $\Delta=3.12$ , $\sqrt{\Delta}=1.7663521733\cong 1.77$, and its roots are 
$E_{T_{1}}=0.4913043478\cong 0.49  \   and    \   E_{T_{2}}=2.0304347826\cong 2.03 $. \\
It should be noted that $E_{T_{2}}$ is approximately the sum of the levels of effort  of the players in Nash equilibrium and these two roots indicate that the amount of harvest in the Nash equilibrium  can be obtained with  $E_{T_{2}}$ and with less effort $E_{T_{1}}$.
Since $E_{T_{2}}-E_{T_{1}}\cong 1.54$, then in $E_{T_{2}}$ that  players are  doing  a static game with complete information 
thay are wasting   $1.54$ units of effort because the mass of fish is reduced and fishing is hardly done
and players should spend more effort to  harvest more. 
Therefore, by considering certain values as parameters, in order to make a clearer analysis, we conclude when
 the harvest  function is  nonlinear, doing a static game causes a waste of energy and money of both players.
Therefore, doing  a cooperative game with the agreement of the parties  will be for the benefit of each player.


\end{document}